\documentclass[11pt]{article}

\usepackage{amsmath,amssymb,latexsym,amsfonts,enumerate,theorem,longtable,amscd}

\newtheorem{theorem}{Theorem}[section]

\newtheorem{lemma}[theorem]{Lemma}
\newtheorem{defn}[theorem]{Definition}
{

\newtheorem{remark}[theorem]{Remark}

\newtheorem{prop}[theorem]{Proposition}

} 
\numberwithin{equation}{section}

\newcommand{\C}{\mathbb C}

\newcommand{\QED}{\ifhmode\unskip\nobreak\fi\quad\ensuremath\square}

\allowdisplaybreaks

\title{{A new construction of $\tilde{D}_5$-singularities and generalization of Slodowy slices \\ 
       } 
\author{Kazunori Nakamoto\footnote{Center for Life Science Research, University of Yamanashi, Yamanashi 409-3898, 
Japan, e-mail: nakamoto@yamanashi.ac.jp} 
 and Meral Tosun\footnote{Galatasaray University, Istanbul, Turkey, e-mail: mtosun@gsu.edu.tr}}
\date{}
}

\begin{document}

\maketitle

\begin{abstract}
Any simple elliptic singularity of type $\tilde{D}_5$ can be obtained by taking the 
intersection of the nilpotent variety and the $4$-dimensional "good slices" in the semi-simple Lie algebra 
${\frak sl}(2, {\mathbb C}) \oplus {\frak sl}(2, {\mathbb C})$.  We describe these new slices  purely by  the structure of the Lie algebra. 
We also construct the semi-universal deformation spaces of $\tilde{D}_5$-singularities
by using the 4-dimensional "good slices". 
\end{abstract}

\footnote[0]{The first author is partially supported by Grant-in-Aid for Young Scientists (B) (No. 19740007)}
\footnote[0]{The second author is supported by TUBITAK 1001 projects with the number 109T667. }
\footnote[0]{Key words: Simple elliptic singularities, $\tilde{D}_5$-singularities, 
Lie algebras, Deformations, Bj\"orner-Welker sequence.}

\section{Introduction}

The first relation established between singularities and Lie algebras is due to the work of Du Val in  \cite{DV} and Artin in  \cite{Artin}; they show that the Dynkin diagrams appearing as the root systems of the semisimple Lie algebras of type ADE are the dual graphs of the minimal resoluion of the isolated complex surface singularities. 
 
The second relation between these two objects is the Grothendieck-Brieskorn theory (see \cite{grothen}, \cite{Brieskorn1}). As our work is based on that theory let us recall it briefly: Let ${\frak g}$ be a simple Lie algebra over $\C$ and $G$ the adjoint group of ${\frak g}$.  The classical theorem of Chevalley states the isomorphism 
$\C[{\frak g}]^G\cong \C[{\frak h}]^ W$ where ${\frak h}$ is a Cartan subalgebra and $W$ is the Weyl group  of ${\frak g}$ (\cite{CH}).  This induces the morphism  
 $\gamma : {\frak g}\rightarrow {\frak h}/W$ called the  adjoint quotient  map which is the key point of the Grothendieck-Brieskorn theory.   When ${\frak g}$ is of ADE type, the adjoint orbit of codimension two in the nilpotent variety $\gamma^{-1}(0) $ of ${\frak g}$  has a rational double singularity  (\cite{Brieskorn1}). 
In \cite{Slodowy1}, Slodowy proved this fact in a different way and introduced an important tool, Slodowy slices, by which we can define the singularities of the adjoint orbits in ${\frak g}$  purely by using the corresponding Lie algebra: For a Jacobson-Morozov $sl_2$-triple $\{x,y,h\}$, a Slodowy slice to the adjoint orbit of a nonzero nilpotent element $x$ in ${\frak g}$ is the affine space $S=x+Ker(ad y)$.  This exceptional  choice of slices are used in many works in mathematics (see \cite{P}, \cite{GHC}). 

Here using the idea of the Grothendieck-Brieskorn theory we relate a new Lie algebra  with a different class of singularities of complex surfaces. 
This new Lie algebra is  ${\frak sl}(2, {\mathbb C}) \oplus {\frak sl}(2, {\mathbb C})$ and corresponds to the simple elliptic singularities  of type $\tilde{D}_5$. 
More precisely the simple elliptic singularities of surfaces are defined by K. Saito in  \cite{Saito2}, named as $\tilde D_5$,  $\tilde E_6$, $\tilde E_7$ and $\tilde E_8$.  For these singularities, the exceptional divisor in the minimal resolution is a smooth elliptic curve of self-intersection respectively $-4$, $-3$, $-2$ and $-1$. The $\tilde{D}_5$-singularities are defined by two quadrics in $\C^4$. There have been several attempts to construct a Lie algebra corresponding to the simple elliptic singularities of types  $\tilde E_6$, $\tilde E_7$ and $\tilde E_8$ and their deformations  (see  \cite{Saito2},  \cite{Looijenga2},   \cite{Looijenga3}, \cite{Pink}, \cite{Mer}). First, Saito constructed the extended affine and elliptic root systems and then many mathematicians tried to construct a Lie algebra admitting a given root system as its real roots (\cite{Slodowy2}).   We can send the interested reader to the excellent introduction of  \cite{Saito3} for the historical developement of the subject. 

Here we deal with the $\tilde D_5$-singularities and our construction contrasts with the one of Helmke and Slodowy \cite{HS} who used an infinite dimensional loop group. As in the case of the simple singularities of surfaces, a $\tilde{D}_5$-singularity can be  obtained by taking the intersection of the nilpotent variety $\gamma^{-1}(0)$ of 
${\frak g}:={\frak sl}(2, {\mathbb C}) \oplus {\frak sl}(2, {\mathbb C})$ and a $4$-dimensional "good" slice passing through the origin in the Lie algebra 
${\frak g}$ (\cite{CRAS}).  For choosing "good" slices, we introduce the notion of $2$-dimensional "good" subspace of ${\frak sl}(2, {\mathbb C}) \oplus {\frak sl}(2, {\mathbb C})$ which is spanned by $2$ vectors $x, y$ of special forms.  The normalization of the good subspaces gives us  a simple coordinate $(p, q)$ for $2$-dimensional good subspaces. Taking the $4$-dimensional subspace orthogonal to a $2$-dimensional good subspace with respect to the Killing form, 
we obtain a good slice. 

Moreover the $j$-invariant of the exceptional curve in the 
minimal resolution of the $\tilde{D}_5$-singularity can be calculated with respect to $(p, q)$. 
The surjectivity of the $j$-function leads us to construct any $\tilde{D}_5$-singularity by our method. 
Then we compute explicitly semi-universal deformations of $\tilde{D}_5$-singularities by deforming the adjoint quotient map 
and the "good" slices of ${\frak sl}(2, {\mathbb C}) \oplus {\frak sl}(2, {\mathbb C})$. 
This fact implies that the deformations of the adjoint quotient map and the slices 
have enough information to induce deformations of $\tilde{D}_5$-singularities.  
We give a negative answer in Remark  (\ref{remark:bjorner}) to  the natural question whether the Lie algebra ${\frak g}:={\frak sl}(2, {\mathbb C}) \oplus {\frak sl}(2, {\mathbb C})\oplus \ldots \oplus {\frak sl}(2, {\mathbb C})$ ($n$ times with $n\geq 3$) gives semi-universal deformations of the surface singularity obtained by the same method. 
 In Appendix, we discuss quadratic forms and $\tilde{D}_5$-singularities. We calculate the $j$-invariant of the exceptional curve in the minimal resolution of a $\tilde{D}_5$-singularity.  

Throughout this article, our notation will be ${\frak g} := {\frak sl}(2, {\mathbb C}) \oplus 
{\frak sl}(2, {\mathbb C})$ and $G := {\rm SL}(2, {\mathbb C})\times 
{\rm SL}(2, {\mathbb C})$.  Recall that the group $G$ acts canonically on the Lie algebra ${\frak g}$ of $G$.

\section{Construction of $\tilde{D}_5$-singularities} 

To construct $\tilde{D}_5$-singularities 
by using the Lie algebra ${\frak g} = {\frak sl}(2, {\mathbb C}) \oplus 
{\frak sl}(2, {\mathbb C})$, we first introduce the followings:

\begin{defn}\rm
Let  $V \subset {\frak g}$ be a $2$-dimensional subspace. We say that $V$ is a {\it good}  subspace if 
for a basis $x, y \in V$, we have 
$x = (x_s, x_n)$ and $y = (y_n, y_s)$, where 
$x_s, y_s \in {\frak sl}(2, {\mathbb C})$ are non-zero 
semi-simple elements and $x_n, y_n \in 
{\frak sl}(2, {\mathbb C})$ are non-zero nilpotent elements.  
\end{defn}

\begin{defn}\rm 
Let ${\mathcal S} \subset {\frak g}$ be a $4$-dimensional subspace. We say that 
${\mathcal S}$ is a {\it good slice} if there exists a $2$-dimensional good 
subspace $V \subset {\frak g}$ such that ${\mathcal S} = V^{\perp} = 
\{ z \in {\frak g} \mid \langle z, v \rangle =0 \mbox{ for each } v \in V \}$, 
where $\langle \ast, \ast \rangle$ is the Killing form of ${\frak g}$. 
Note that $\langle (X_1, X_2), (Y_1, Y_2) \rangle = 4 ({\rm tr}(X_1Y_1) + {\rm tr}(X_2Y_2))$ 
for $(X_1, X_2), (Y_1, Y_2) \in {\frak g}$. 
\end{defn}

\begin{lemma}\label{lemma:normalization}  
Let $V \subset {\frak g}$ be a $2$-dimensional good subspace. 
Then there exists $g = (P, Q) \in G$ such that 
$g^{-1} V g$ has a basis $x = (x_s, x_n)$ and 
$y = (y_n, y_s)$ with the following properties: 
\begin{enumerate} 
\item[(1)]   
$x_s = 
\left( 
\begin{array}{cc}
1 & 0 \\
0 & -1 \\
\end{array}
\right), 
y_s = 
\left( 
\begin{array}{cc}
1 & 0 \\
0 & -1 \\
\end{array}
\right)$.   \\
\item[(2)] 
$x_n = 
\left( 
\begin{array}{cc}
p & 1 \\
-p^2 & -p \\
\end{array}
\right) 
\mbox{  for some } 
p \in {\mathbb C}, 
\mbox{ or }  
x_n = 
\left( 
\begin{array}{cc}
0 & 0 \\
1 & 0 \\
\end{array}
\right)$.   \\
\item[(3)] 
$y_n = 
\left( 
\begin{array}{cc}
q & 1 \\
-q^2 & -q \\
\end{array}
\right) 
\mbox{  for some } 
q \in {\mathbb C}, 
\mbox{ or }  
y_n = 
\left( 
\begin{array}{cc}
0 & 0 \\
1 & 0 \\
\end{array}
\right)$.   
\end{enumerate} 
\end{lemma}

\noindent {\it Proof}.
Since $V$ is a good subspace in ${\frak g}$, we can choose a basis 
$x = (x_s, x_n)$ and $y = (y_n, y_s)$ of $V$, where 
$x_s, y_s \in {\frak sl}(2, {\mathbb C})$ are non-zero 
semi-simple elements and $x_n, y_n \in 
{\frak sl}(2, {\mathbb C})$ are non-zero nilpotent elements.  
Remplacing $x, y$ by $ax, by$ with suitable $a, b \in 
{\mathbb C}^{\times}$ respectively, we may assume that 
$\det x_s = \det y_s = -1$. Hence we can  
take $P', Q' \in {\rm SL}(2, {\mathbb C})$ such that 
$P'^{-1}x_s P' = 
\left( 
\begin{array}{cc}
1 & 0 \\
0 & -1 \\
\end{array}
\right)$ and  
$Q'^{-1}y_s Q' = 
\left( 
\begin{array}{cc}
1 & 0 \\
0 & -1 \\
\end{array}
\right)$. 

Let $x_n  = 
\left( 
\begin{array}{cc}
\alpha & \beta \\
\gamma & -\alpha \\
\end{array}
\right)$. 
First, assume that $\beta \neq 0$. 
Then putting $Q = Q' 
\left( 
\begin{array}{cc}
\sqrt{\beta} & 0 \\
0 & 1/\sqrt{\beta} \\
\end{array}
\right)$, we have 
$Q^{-1} y_s Q = 
\left( 
\begin{array}{cc}
1 & 0 \\
0 & -1 \\
\end{array}
\right)$ and 
$Q^{-1} x_n Q = 
\left( 
\begin{array}{cc}
\alpha & 1 \\
\beta\gamma & -\alpha \\
\end{array}
\right)$. 
Since $\det x_n =0$, 
we can put $\alpha = p$ and $\beta\gamma = -p^2$. 

Next, assume that $\beta = 0$. 
Since $\det x_n = 0$ and $x_n \neq 0$,  
by putting $Q =  Q' 
\left( 
\begin{array}{cc}
1/\sqrt{\gamma} & 0 \\
0 & \sqrt{\gamma} \\
\end{array}
\right)$, we have 
$Q^{-1} y_s Q= 
\left( 
\begin{array}{cc}
1 & 0 \\
0 & -1 \\
\end{array}
\right)$ and 
$Q^{-1} x_n Q = 
\left( 
\begin{array}{cc}
0 & 0 \\
1 & 0 \\
\end{array}
\right)$. 

In a similar way, let us take $P$ by modifying $P'$ 
such that 
$P^{-1} x_s P= 
\left( 
\begin{array}{cc}
1 & 0 \\
0 & -1 \\
\end{array}
\right)$ and 
$P^{-1} y_n P = 
\left( 
\begin{array}{cc}
q & 1 \\
-q^2 & q \\
\end{array}
\right)$ for some $q \in {\mathbb C}$, 
or $P^{-1} y_n P = \left( 
\begin{array}{cc}
0 & 0 \\
1 & 0 \\
\end{array}
\right)$. 
This completes the proof. 
\hfill $\Box$ 

\bigskip

Then for a $2$-dimensional good subspace $V$ of ${\frak g}$, it is enough to consider only the case 
$V = {\mathbb C}x + {\mathbb C}y$ with the properties $(1), (2)$ and $(3)$ 
in Proposition \ref{lemma:normalization}. 

\bigskip 

Let us denote by ${\mathcal N}({\frak g})$  the nilpotent variety of ${\frak g}$ which is 
 described as 
\[
{\mathcal N}({\frak g}) = \left\{ 
\begin{array}{c|c}
\left( \left(
\begin{array}{cc}
a & b \\
c & -a \\
\end{array}
\right), 
\left( 
\begin{array}{cc}
d & e \\
f & -d \\
\end{array}
\right)
 \right) \in {\frak g} &
a^2+bc=d^2+ef=0 
\end{array}
\right\}.
\]
For a $2$-dimensional good subspace 
$V_{(x, y)} = {\mathbb C}x + {\mathbb C}y \subset {\frak g}$, 
we can define the good slice ${\mathcal S}_{(x, y)} = 
V_{(x, y)}^{\perp}$. 
Then we define $X_{(x, y)} := {\mathcal N}({\frak g}) \cap {\mathcal S}_{(x, y)}$. 
The germ $(X_{(x, y)}, 0)$ is a $\tilde{D}_5$-singularity for general 
$x=(x_s, x_n)$ and $y=(y_n, y_s)$. 

\begin{remark}\label{remark:1}\rm
Let $x=(x_s, x_n)$ and $y=(y_n, y_s)$ with 
$x_s = y_s = \left( 
\begin{array}{cc}
1 & 0 \\
0 & -1 \\
\end{array}
\right)$,  
$x_n = 
\left( 
\begin{array}{cc}
p & 1 \\
-p^2 & -p \\
\end{array}
\right)$, and 
$y_n = 
\left( 
\begin{array}{cc}
q & 1 \\
-q^2 & -q \\
\end{array}
\right)$ for $p, q \in {\mathbb C}$.   

For $Z := \left( \left( 
\begin{array}{cc}
a & b \\
c & -a \\
\end{array}
\right), \left(
\begin{array}{cc}
d & e \\
f & -d \\
\end{array}
\right) \right) 
 \in {\frak g}$, 
Z is contained in ${\mathcal S}_{(x, y)}$ if and only if 
$2a+2pd-p^2e+f=0$ and $2qa-q^2b+c+2d=0$. 

Furthermore, the germ $X_{(x, y)}$ is described as 
\[
\begin{array}{ccl}
X_{(x, y)} & = & \left\{ (a, b, c, d, e, f) \in {\mathbb C}^6 \; 
\begin{array}{|l} 
 2a+2pd-p^2e+f=0 \\
 2qa-q^2b+c+2d=0 \\
 a^2+bc= d^2+ef=0 \\
\end{array} 
\right\} \smallskip \\
 & = &  
\left\{ (a, b, d, e) \in {\mathbb C}^4 \; 
\begin{array}{|l} 
 a^2-2qab+q^2b^2-2bd=0 \\
 -2ae+d^2-2pde+p^2e^2=0 \\
\end{array} 
\right\}. \\
\end{array}
\]
\end{remark}

\begin{prop}\label{prop:jcalculated}
Let $x=(x_s, x_n)$ and $y=(y_n, y_s)$ be as in Remark \ref{remark:1}. 
Set $t=pq$. The germ $(X_{(x, y)}, 0)$ is a 
$\tilde{D}_5$-singularity for $t \neq 0, 1/4$. 
For the minimal resolution $\mu : \tilde{X}_{(x, y)} \to X_{(x, y)}$, the $j$-invariant of the exceptional curve
$E_{(x, y)}$  of $\mu$ is given by 
\[\displaystyle 
j(E_{(x, y)}) = \frac{256(t^6-12t^5+51t^4-88t^3+51t^2-12t+1)}{t^4(1-4t)}.  
\] 
\end{prop}

{\it Proof}. 
Let \[
X = 
\left(
\begin{array}{cccc}
1 & -q & 0 & 0 \\
-q & q^2 & -1 & 0 \\
0 & -1 & 0 & 0 \\
0 & 0 & 0 & 0 \\
\end{array} 
\right) \text{ and }  
Y = 
\left(
\begin{array}{cccc}
0 & 0 & 0 & -1 \\
0 & 0 & 0 & 0 \\
0 & 0 & 1 & -p \\
-1 & 0 & -p & p^2 \\
\end{array} 
\right). 
\]
Then $X_{(x, y)} \cong S_{(X, Y)} := 
\{ v={}^t (x, y, z, w) \in {\mathbb C}^4 \mid {}^t v X v = {}^t v Y v = 0 \}$. 

For each $p, q \in {\mathbb C}$ there exists  
$\left( 
\begin{array}{cc}
a & b \\
c & d \\
\end{array}
\right) \in {\rm GL}_2({\mathbb C})$ 
such that $aX+bY \in {\rm GL}_4({\mathbb C})$.  
Let us choose $P \in {\rm GL}_4({\mathbb C})$ such that 
${}^{t}\!P(aX+bY)P = I_4$. 
Putting $Y' := {}^{t}\!P(cX+dY)P$, we have 
$S_{(X, Y)} \cong S_{(aX+bY, cX+dY)} \cong 
S_{(I_4, Y')}$. 
Since $Y' = P^{-1} (aX+bY)^{-1} (cX+dY) P$, the germ 
$(X_{(x, y)}, 0)$ is $\tilde{D}_5$-singularity if and only if 
the characteristic polynomial of $(aX+bY)^{-1} (cX+dY)$  
has no multiple root by Lemma \ref{lemma:mulroot}. 
Using Theorem \ref{th:jinvariant}, we can calculate  
the $j$-invariant $j(Y')$ of $E_{(x, y)}$     
in terms of the coefficients of 
the characteristic polynomial of $(aX+bY)^{-1} (cX+dY)$. 
\hfill $\Box$

\begin{remark}\label{remark:jcalculated}\rm
If $t=0$ or $t=1/4$, then $(X_{(x, y)}, 0)$ is not a $\tilde{D}_5$-singularity. 
In the case that 
$x_n = \left( \begin{array}{cc} 
0 & 0 \\
1 & 0 \\
\end{array} 
\right)$ or 
$y_n= \left( \begin{array}{cc} 
0 & 0 \\
1 & 0 \\
\end{array} 
\right)$, then $(X_{(x, y)}, 0)$ is not a $\tilde{D}_5$-singularity  
either. 
\end{remark}

\begin{remark}\rm
Any $\tilde{D}_5$-singularity can be obtained by 
choosing suitable $x, y$ as in Remark \ref{remark:1} 
since
\[
\begin{array}{ccccc}
j  & : & {\mathbb C}\times {\mathbb C} & \to & 
{\mathbb C}\cup \{ \infty \} \\
 & & (p, q) & \mapsto & j(E_{(x, y)})
\end{array} 
\]
is surjective. 
\end{remark}

\section{Semi-universal deformations}

Now we want to construct a semi-universal deformation of 
the ${\tilde D}_5$-singularity $(X_{(x, y)},0)$ by using 
the Lie algebra ${\frak g} = {\frak sl}(2, {\mathbb C}) \oplus 
{\frak sl}(2, {\mathbb C})$ in terms of Lie algebras. 

\bigskip  

Recall that the germ $(X_{(x,y)}, 0)$ is described as 
\[
X_{(x, y)} =   
\left\{ (a, b, d, e) \in {\mathbb C}^4 \; 
\begin{array}{|l} 
 g_1 := a^2-2qab+q^2b^2-2bd=0 \\
 g_2 :=-2ae+d^2-2pde+p^2e^2=0 \\
\end{array} 
\right\}.
\]
Let ${\mathcal O}_{(X_{(x,y)}, 0)}$ be the local ring of the germ $(X_{(x,y)}, 0)$.
For obtaining a semi-universal deformation space of 
$(X_{(x,y)}, 0)$, we calculate the vector space 
$T^1 = {\mathcal O}_{(X_{(x,y)}, 0)}^2/M'$, where 
$M'$ is the ${\mathcal O}_{(X_{(x,y)}, 0)}$-submodule 
of ${\mathcal O}_{(X_{(x,y)}, 0)}^2$ generated by 
the $4$ vectors: 
$\displaystyle (\frac{\partial g_1}{\partial a},
\frac{\partial g_2}{\partial a})$, $\displaystyle (\frac{\partial
g_1}{\partial b}, \frac{\partial g_2}{\partial b})$, $\displaystyle
(\frac{\partial g_1}{\partial d}, \frac{\partial g_2}{\partial d})$
and $\displaystyle (\frac{\partial g_1}{\partial e}, \frac{\partial
g_2}{\partial e})$. (For details, see \cite{Tjurina}.) 
Let ${\mathcal O} = {\mathbb C}\{a, b, d, e\}$. 
Since ${\mathcal O}_{(X_{(x,y)}, 0)} = {\mathcal O}/(g_1, g_2)$, 
we have $T^1 = {\mathcal O}^2/M$, where $M$ is the ${\mathcal O}$-submodule 
of ${\mathcal O}^2$ generated by the following vectors: 
\[
v_1  =  (g_1, 0), \;  v_2 = (g_2, 0), \; v_3 = (0, g_1), \;  v_4=(0, g_2), 
\]
\begin{eqnarray*} 
v_5 & = &  \frac{1}{2}(\frac{\partial g_1}{\partial a}, \frac{\partial g_2}{\partial a}) 
= (a-qb, -e),  \\  
v_6 & = &  -\frac{1}{2}(\frac{\partial g_1}{\partial b}, \frac{\partial g_2}{\partial b}) 
= (qa-q^2b+d, 0), \\ 
v_7 & = & -\frac{1}{2}(\frac{\partial g_1}{\partial d}, \frac{\partial g_2}{\partial d}) 
= (b, -d+pe), \\ 
v_8 & = & -\frac{1}{2}(\frac{\partial g_1}{\partial e}, \frac{\partial g_2}{\partial e})
= (0, a+pd-p^2e). \\ 
\end{eqnarray*} 
If $t=pq$ is equal to $0$ or $1/4$, then $(X_{(x,y)}, 0)$ is not 
a $\tilde{D}_5$-singularity by Remark \ref{remark:jcalculated}. 
Hence we assume that $t \neq 0, 1/4$.

\begin{lemma}\label{lemma:basis} 
If $t \neq 0, 1/4$, then 
the vector space $T^{1}$ has a basis 
\[
(1, 0), (b, 0), (e, 0), (0, 1), (0, b), (0, ae), (0, e).
\]
\end{lemma}

\noindent {\it Proof}. 
In the following, it is useful to consider 
the lexicographic order on ${\mathcal O}$ with $a > b > d > e$ and 
the monomial order $>_{POT}$ on ${\mathcal O}^2$ such that $(1, 0) >_{POT} (0, 1)$ 
(cf. \cite[Chapter 5]{CLO}). 
It is easy to check the following equalities. 
\begin{eqnarray*}
v_9 & := & v_5+qv_7 = (a, -e-qd+pqe), \\
v_{7} & = & (b, -d+pe), \\ 
v_{10} & := & v_6-qv_5 = (d, qe), \\
v_{11} & := & v_2+2ev_9 +(2pe-d)v_{10} = (p^2e^2, -3qde+(4pq-2)e^2), \\
v_{8} & = & (0, a+pd-p^2e), \\
v_{12} & := & v_3+(-a+2qb+pd-p^2e)v_8 \\ 
 & = & (0, q^2b^2+(2pq-2) bd - 2p^2q be +p^2d^2 -2p^3 de + p^4 e^2), \\
v_{13} & := & av_7+(d-pe)v_8 -bv_9  \\ 
& = & (0, qbd + (1-pq)be +pd^2-2p^2de+p^3e^2), \\ 
v_{14} & := & - dv_{7} + bv_{10} = (0, qbe+d^2-pde), \\ 
v_{15} & := & v_{4} + 2ev_{8} = (0, d^2-p^2e^2), \\
v_{16} & := & qv_4 + 3qev_8 + dv_9 - av_{10} \\
& = & (0, -2p^2q e^2 + (2pq-1)de), \\
v_{17} & := & \frac{1}{4pq-1} \{ 4p^2q^2 e v_{15} -(2pq-1)d v_{16} - 2 p^2 q ev_{16} \} \\ 
 & = & (0, d^2e), \\
v_{18} & := & \frac{1}{2p^2q} (-dv_{16}+(2pq-1)v_{17}) \\
 & = & (0, de^2). 
\end{eqnarray*}
All $v_i$ ($1 \le i \le 18$) are contained in $M$. 
We see that the following $7$ vectors span $T^1 = {\mathcal O}^2/M$: 
\[
(1, 0), (e, 0), (0, 1), (0, b), (0, d), (0, e), (0, de).
\] 
Since $\dim T^1 =7$ for $\tilde{D}_5$-singularities, 
the vectors above are a basis of $T^1$. By the equalities  
\begin{eqnarray*}
(0, d) & = & -v_7+(b, 0)+p(0, e) \;\text{ and}  \\
(0, de) & = & 2qe v_{8} -2q (0, ae) -v_{16},    
\end{eqnarray*} 
we have another basis of $T^1$: 
\[
(1, 0), (b, 0), (e, 0), (0, 1), (0, b), (0, ae), (0, e).
\]
This completes the proof. 
\hfill $\Box$ 


\bigskip 

To construct a semi-universal deformation of the $\tilde{D}_5$-singularity 
$(X_{(x, y)}, 0)$ purely by the Lie algebra point of view  
${\frak g} = {\frak sl}(2, {\mathbb C}) \oplus 
{\frak sl}(2, {\mathbb C})$, consider the adjoint quotient of ${\frak g}$ is described as 
\[
\begin{array}{ccccl}
f & : & {\frak g} & \to & {\frak g}/\!/G \cong {\frak h}/W \cong {\mathbb C}^2 \\
  &   &  z = (z_1, z_2) & \mapsto & (\det z_1, \det z_2).  \\
\end{array}
\]
where ${\frak h} := \left\{ 
\left( \left( 
\begin{array}{cc}
a & 0 \\
0 & -a \\
\end{array}
\right), \left( 
\begin{array}{cc}
d & 0 \\
0 & -d \\
\end{array} 
\right) \right) \in {\frak g} 
\right\}$ be a Cartan subalgebra of ${\frak g}$. 
and $W$ is the Weyl group of ${\frak g}$ isomorphic to ${\mathbb Z}/2{\mathbb Z} \times 
{\mathbb Z}/2{\mathbb Z}$.

\noindent Set $x_{\infty} := 
\left( \begin{array}{cc} 
0 & 0 \\
1 & 0 \\
\end{array} 
\right)$ and 
$y_{\infty} := 
\left( \begin{array}{cc} 
0 & 0 \\
1 & 0 \\
\end{array} 
\right)$.  
\noindent Let us deform the adjoint quotient $f$ by $(\alpha, \beta) \in {\mathbb C}^2$ 
as follows: 
\[
\begin{array}{ccccl}
f_{(\alpha, \beta)} & : & {\frak g} & \to & {\frak g}/\!/G \cong {\frak h}/W \cong {\mathbb C}^2 \\
  &   &  z=(z_1, z_2) & \mapsto & (\det z_1 + \alpha \langle z, (0, x_{\infty})\rangle,  
\det z_2 + \beta \langle z, (y_{\infty}, 0) \rangle).   \\
\end{array}
\]

Next, recall the slice 
${\mathcal S}_{(x, y)} = \{ z \in {\frak g} \mid \langle z, x \rangle 
= \langle z, y \rangle = 0 \}$ and deform it by $(\gamma, \delta, \varepsilon) \in 
{\mathbb C}^3$ as follows: 
\[
{\mathcal S}_{(x, y)}(\gamma, \delta, \varepsilon) := 
\{ z \in {\frak g} \mid \langle z, x \rangle + \gamma \langle z, (x_s, 0) \rangle = 4\delta, 
\; \langle z, y \rangle = 4\varepsilon  
\}. 
\] 

\noindent We then consider $S := {\mathbb C}^2 \times {\mathbb C}^3 \times {\frak h}/W 
\cong {\mathbb C}^7$ as the base space and  
\[
{\mathcal X} := \{ (Z, \alpha, \beta, \gamma, \delta, \varepsilon, 
\lambda, \mu) \in {\frak g}\times S \mid f_{(\alpha, \beta)}(Z)=(\lambda, \mu), 
Z \in {\mathcal S}_{(x, y)}(\gamma, \delta, \varepsilon) 
\} 
\]
as the total space. We denote by $\pi: {\mathcal X} \to S$ the second projection. 

\begin{theorem}
The projection $\pi : ({\mathcal X}, 0) \to (S, 0)$ gives a 
semi-universal deformation of $(X_{(x, y)}, 0)$ for 
$pq \neq 0, 1/4$. 
\end{theorem}

\noindent {\it Proof}. 
The deformed adjoint quotient $f_{(\alpha, \beta)}$ 
is described as 
\[
f_{(\alpha, \beta)}(z) = (-a^2-bc+4\alpha e, -d^2-ef+4\beta b). 
\]
The deformed slice ${\mathcal S}_{(x, y)}(\gamma, \delta, \varepsilon)$ is described as 
\[
{\mathcal S}_{(x, y)}(\gamma, \delta, \varepsilon) = \{ z \in {\frak g} \mid 
2a+2\gamma a+ 2pd-p^2e+f = \delta, 2qa-q^2b+2d+c=\varepsilon \}. 
\]
Hence the total space ${\mathcal X}$ is 
\[
{\mathcal X} = \left\{ (a, b, d, e) \in {\mathbb C}^4 \;  
\begin{array}{|l}
(g_1, g_2) 
- 4\alpha (e, 0)   
- 4\beta (0, b)  
-2 \gamma (0, ae) \\   
+\delta (0, e) 
+\varepsilon (b, 0) 
+\lambda (1, 0) 
+\mu (0, 1) 
= (0, 0) \\
\end{array}
\right\}. 
\]
Since the vectors 
\[
(1, 0), (b, 0), (e, 0), (0, 1), (0, b), (0, ae), (0, e).
\]
are a basis of $T^1$ by Lemma \ref{lemma:basis}, the projection 
$\pi : ({\mathcal X}, 0) \to (S, 0)$ gives a 
semi-universal deformation. 
\hfill $\Box$

\begin{remark}\rm
Our construction of $\tilde{D}_5$-singularities and their 
semi-universal deformation spaces is given in terms of Lie algebras. 
However, it seems that the construction of semi-universal deformation spaces 
is a little bit artificial, because we can not answer why 
$\dim T^1 = 7$ essentially.  
It is more natural to say that we can obtain versal deformations of 
$\tilde{D}_5$-singularities only by deforming the adjoint quotients and 
slices.  
\end{remark} 

\begin{remark}\label{remark:bjorner}
\rm
Let ${\frak g}_m := \overbrace{{\frak sl}(2, {\mathbb C})\oplus
{\frak sl}(2, {\mathbb C})\oplus \cdots \oplus {\frak sl}(2, {\mathbb
C})}^{m}$. Let us consider a general slice ${\mathcal S}$ at $0$ 
such that $({\mathcal N}({\frak g}_m) \cap {\mathcal S}, 0)$ is an  
isolated surface singularity.  
Here we denote by ${\mathcal N}({\frak g}_m)$ the nilpotent variety 
of ${\frak g}_m$. We have
\begin{eqnarray*}
BW(m):=dim T^1({\mathcal N}({\frak g}_m) \cap {\mathcal S})  & = & \sum_{i=3}^{m+2} 
\binom{m+2}{i}\binom{i-1}{2} \\ 
& = & 2^{m-1}(m^2-m+2)-1,  
\end{eqnarray*} 
although we don't prove here. 
The sequence $BW(m)$ is called Bj\"orner-Welker sequence, which 
is equal to ${\rm rank} H^1 (M_{m+2, 3}^{\mathbb R})$, where  
$M_{m+2, 3}^{\mathbb R}$ is the complement of  
\[
V_{m+2, 3}^{\mathbb R} = \{ (x_i) \in {\mathbb R}^{m+5} \mid x_{i_1} = x_{i_2} = x_{i_3} 
\text{ for some indices } i_1, i_2, i_3 \}
\] 
in ${\mathbb R}^{m+5}$ (See \cite{BW}). 

Note that $\dim {\mathcal S}=m+2$. Let ${\rm Aff}({\frak g}_m, m+2)$ be 
the variety of $(m+2)$-dimensional affine subspaces of ${\frak g}_m$. 
The variety ${\rm Aff}({\frak g}_m, m+2)$ can be embedded in 
the Grassmann variety ${\rm Grass}(\dim {\frak g}_m + 1, m+2)$. 
Since $\dim {\rm Aff}({\frak g}_m, m+2) = \dim 
{\rm Grass}(\dim {\frak g}_m + 1, m+2) = (2m-1)(m+2)$, 
the dimension of deformations of slices has polynomial order. 
We can imagine that the dimension of deformations of the adjoint quotient 
has polynomial order, too. On the other hand, $BW(m)$ has exponential order.  
Hence, the deformations of the adjoint quotient and slices have no potential to 
induce the deformation of the singularity ${\mathcal N}({\frak g}_m) \cap {\mathcal S}$.  
In general, we can not construct a semi-universal deformation space of 
${\mathcal N}({\frak g}_m) \cap {\mathcal S}$ by our method. 
\end{remark} 

\section{Appendix}

In this appendix, we discuss quadratic forms and $\tilde{D}_5$-singularities. 
We calculate the $j$-invariant of the exceptional curve in the minimal 
resolution of a $\tilde{D}_5$-singularity. 

First, we classify the singularities defined by
two quadratic equations in ${\mathbb C}^4$ written by
\[
f(x, y, z, w)={}^t v X v, \quad  g(x, y, z, w)= {}^t vYv,
\]
where $X$ and $Y$ are $4 \times 4$ symmetric matrices and $v = {}^t (x,
y, z, w)$. 
Set $S_{(X, Y)} := \{ (x, y, z, w)
\in {\mathbb C}^4 \mid f(x, y, z, w)=g(x, y, z, w)=0 \}$.

\begin{defn}\rm
Let ${\rm Sym}_n({\mathbb C})$ be the set of $n \times n$ symmetric matrices over ${\mathbb C}$.
For pairs $(X, Y), (X', Y') \in {\rm Sym}_n({\mathbb C})\times {\rm Sym}_n({\mathbb C})$
we say that $(X, Y) \sim (X', Y')$ if there exists $P \in {\rm GL}_n({\mathbb C})$ such that
${}^{t}\!PXP = X'$ and ${}^{t}\!PYP=Y'$.
\end{defn}

\begin{lemma}\label{lemma:basicsymmetric}
Let $X \in {\rm Sym}_n({\mathbb C})$.
If $X \neq 0$, then there exists $v_0 \in {\mathbb C}^n$
such that ${}^t v_0 X v_0 = 1$.
\end{lemma}

\noindent {\it Proof}. By the assumption, there exist $v, w \in {\mathbb
C}^n$ such that ${}^t v X w \neq 0$. If ${}^t v X v \neq 0$ or ${}^t
w X w \neq 0$, then we can obtain $v_0$ by a suitable scalar
multiplication on $v$ or $w$. Suppose that ${}^t v X v = 0$ and
${}^t w X w = 0$. Then ${}^t (v+w) X (v+w) = 2 {}^t vXw \neq 0$. In
this case we can get $v_0$ by a suitable scalar 
multiplication on $v+w$. 
\hfill $\Box$

\begin{lemma}\label{lemma:normalformsym}
Let $X \in {\rm Sym}_n({\mathbb C})$.
Then there exists $P \in GL(n, {\mathbb C})$ such that
\[
{}^{t}\!PXP =
\left(
\begin{array}{cc}
I_r & 0 \\
0 & 0 \\
\end{array}
\right),
\]
where $r$ is the rank of $X$.
\end{lemma}

\noindent {\it Proof}.
By Lemma \ref{lemma:basicsymmetric} if $X \neq 0$ there exists
$v_1 \in {\mathbb C}^n$ such that ${}^t v_1 X v_1 = 1$.
Let $V_1 := \{ v \in {\mathbb C}^n \mid {}^t v_1 X v = 0 \}$.
By restricting $X$ into $X\!\!\!\mid_{V_1}$, the size of symmetric
matrices decreases. If we continue to apply Lemma \ref{lemma:basicsymmetric}
we get the statement.
\hfill $\Box$

\begin{lemma}\label{lemma:symmat}
Let $X, X' \in {\rm Sym}_n({\mathbb C})$.
Assume that $X$ is diagonalizable.
Then there exists $P \in {\rm GL}_n({\mathbb C})$ such that $P^{-1}XP=X'$
if and only if there exists $Q \in {\rm O}_n({\mathbb C})$ such that
${}^t QXQ=X'$.
Here ${\rm O}_n({\mathbb C}) = \{ Q \in {\rm GL}_n({\mathbb C}) \mid {}^t QQ=I_n \}$.
\end{lemma}

\noindent {\it Proof}.
The "if" part is obvious.
Let us show the "only if" part.
Suppose that there exists $P \in {\rm GL}_n({\mathbb C})$ such that $P^{-1}XP=X'$.
Since $X$ is diagonalizable, we may assume that
$X$ and $X'$ have the same Jordan form
${\rm diag}(\alpha_1, \alpha_2, \ldots, \alpha_n)$.
It suffice to show that there exist $Q, Q' \in {\rm O}_n({\mathbb C})$ such that
${}^t QXQ= {}^t Q'X'Q'= {\rm diag}(\alpha_1, \alpha_2, \ldots, \alpha_n)$.

Let $v_1, v_2, \ldots, v_n \in {\mathbb C}^n$ be eigenvectors of $X$
which belong to eigenvalues $\alpha_1, \alpha_2, \ldots, \alpha_n$,
respectively. If $\alpha_1, \ldots, \alpha_n$ are not distinct, then
we assume that $\alpha_1 = \alpha_2 = \cdots = \alpha_{i_1},
\alpha_{i_{1}+1} = \cdots = \alpha_{i_2}, \ldots, \alpha_{i_{k-1}+1}
= \cdots = \alpha_{n}$. Put $R = ( v_1, v_2, \ldots, v_n )$. Then $R
\in {\rm GL}_n({\mathbb C})$ and $XR=RD$, where $D={\rm
diag}(\alpha_1, \ldots, \alpha_n)$. Multiplying ${}^{t}\!R$ from the
left, we have ${}^{t}\!RXR = {}^{t}\!RRD$. Since ${}^{t}\!RRD = ( {}^t v_i
v_j \alpha_j )_{1 \le i, j \le n}$ and ${}^{t}\!RXR = {}^{t}\!RRD$ is
symmetric, ${}^t v_i v_j \alpha_j = {}^t v_j v_i \alpha_{i}$. If
$\alpha_i \neq \alpha_j$, then ${}^t v_i v_j =0$. Thus we have
\[
{}^{t}\!RR=
\left(
\begin{array}{cccc}
B_1 & 0 & \cdots & 0 \\
0 & B_2 & \cdots & 0 \\
\vdots & \vdots & \ddots & \vdots \\
0 & 0 & \cdots & B_k
\end{array}
\right),
\]
where $B_{\ell}=( {}^t v_{i}v_{j} )_{i_{\ell-1}+1 \le i, j \le
i_{\ell}}$. There exist $T_1, T_2, \ldots T_k$ such that $^{t}
T_{\ell} B_{\ell} T_{\ell}$ are the identity matrices for $1 \le
\ell \le k$. Putting
\[
Q = R
\left(
\begin{array}{cccc}
T_1 & 0 & \cdots & 0 \\
0 & T_2 & \cdots & 0 \\
\vdots & \vdots & \ddots & \vdots \\
0 & 0 & \cdots & T_k
\end{array}
\right),
\]
\noindent we have ${}^t QQ=I_n$ and ${}^t QXQ = {\rm diag}(\alpha_1,
\alpha_2, \ldots, \alpha_n)$. Similarly we can prove the claim for
$X'$. \hfill $\Box$

\begin{prop}
Let $(X, Y), (X', Y') \in {\rm Sym}_n({\mathbb C})\times
{\rm Sym}_n({\mathbb C})$. Assume that $X, X' \in {\rm GL}_n({\mathbb C})$
and $X^{-1}Y$ is diagonalizable.
Then $(X, Y) \sim (X', Y')$ if and only if $X^{-1}Y$ and $X'^{-1}Y'$ have
the same Jordan form.
\end{prop}

\noindent {\it Proof}.
If $(X, Y) \sim (X', Y')$, then there exists $P \in {\rm GL}_n({\mathbb C})$
such that ${}^{t}\!PXP = X'$ and ${}^{t}\!PYP = Y'$.
Since $X'^{-1}Y' = P^{-1}X^{-1}YP$, $X^{-1}Y$ and $X'^{-1}Y'$ have
the same Jordan form.

Conversely, suppose that $X^{-1}Y$ and $X'^{-1}Y'$ have
the same Jordan form. By the assumption that 
$X^{-1}Y$ is diagonalizable, we may denote   
their Jordan normal form by $D = {\rm diag}
(\alpha_1, \ldots, \alpha_n)$.  
Let us take $P \in {\rm GL}_n({\mathbb C})$ such that
${}^{t}\!PXP = I_n$. Then $P = X^{-1}{}^t\!P^{-1}$. 
Since
${}^{t}\!PYP = P^{-1}P{}^{t}\!PYP = P^{-1}(X^{-1}{}^{t}\!P^{-1}){}^{t}\!P YP = P^{-1}X^{-1}YP$,
the symmetric matrix ${}^{t}\!PYP$  has the
Jordan form $D$.
By Lemma \ref{lemma:symmat} there exists $Q \in {\rm O}_n({\mathbb C})$
such that ${}^t Q{}^{t}\!PYPQ \\ = D$.
Hence $(X, Y) \sim (I_n, {}^{t}\!PYP) \sim (I_n, D)$.
Similarly, we have $(X', Y') \sim (I_n, D)$.
Therefore $(X, Y) \sim (X', Y')$.
\hfill $\Box$


\bigskip 

In the sequel, we deal with the case $n=4$.

\begin{defn}\rm
For $(X, Y) \in {\rm Sym}_4({\mathbb C})\times {\rm Sym}_4({\mathbb C})$,
we define $S_{(X, Y)} := \{ (x, y, z, w) \in {\mathbb C}^4 \mid
{}^t v X v = {}^t v Y v = 0, \mbox{ where } v = {}^t (x, y, z, w)
\}$.
\end{defn}

\begin{remark}\rm
Let $(X, Y), (X', Y') \in {\rm Sym}_4({\mathbb C}) \times {\rm Sym}_4({\mathbb C})$.
If $(X, Y) \sim (X', Y')$, then $S_{(X, Y)} \cong S_{(X', Y')}$.
For
$\left( \begin{array}{cc}
a & b \\
c & d \\
\end{array}
\right) \in {\rm GL}_2({\mathbb C})$,
we have $S_{(X, Y)} \cong S_{(aX+bY, cX+dY)}$.
\end{remark}

\begin{lemma}\label{lemma:non-isolated}
If $aX+bY$ is a singular matrix for any $a, b \in {\mathbb C}$,
then $S_{(X, Y)}$ is not an isolated surface singularity.
\end{lemma}

\noindent {\it Proof}.
We only need to investigate the case that $\dim S_{(X, Y)}=2$.
By considering $({}^{t}\!PXP, {}^{t}\!PYP)$ with $P \in {\rm GL}(4, {\mathbb C})$ instead of the
pair $(X, Y)$,
we may assume that
\[
X = \left(
\begin{array}{cccc}
1 & 0 & 0 & 0 \\
0 & 0 & 0 & 0 \\
0 & 0 & 0 & 0 \\
0 & 0 & 0 & 0 \\
\end{array}
\right) \mbox{ or }
\left(
\begin{array}{cccc}
1 & 0 & 0 & 0 \\
0 & 1 & 0 & 0 \\
0 & 0 & 0 & 0 \\
0 & 0 & 0 & 0 \\
\end{array}
\right) \mbox{ or }
\left(
\begin{array}{cccc}
1 & 0 & 0 & 0 \\
0 & 1 & 0 & 0 \\
0 & 0 & 1 & 0 \\
0 & 0 & 0 & 0 \\
\end{array}
\right).
\]
Suppose that $X$ is the first or the second one. There exists
$v_{0} = {}^t (0, 0, z, w) \neq {}^t (0, 0, 0, 0)$ such that ${}^t v_0 Y v_0 = 0$.
For each $\lambda \in {\mathbb C}$, $\lambda v_0$ is a singular point of $S_{(X, Y)}$.
Hence $S_{(X, Y)}$ is not an isolated singularity.

Suppose that $X$ is the third one. Set $Y=(y_{ij})$. By the assumption,
$\det(t X + Y) = y_{44} t^3 + \cdots = 0$ for any $t \in {\mathbb C}$.
So we have $y_{44} =0$.
Set $e_4 = {}^t (0, 0, 0, 1)$. For each $\lambda \in {\mathbb C}$,
${}^t (\lambda e_4) X (\lambda e_4) = {}^t (\lambda e_4) Y (\lambda e_4) = 0$.
Hence $\lambda e_4 \in S_{(X, Y)}$. Furthermore, the Jacobi matrix has rank $<2$
at $\lambda e_4$ for each $\lambda \in {\mathbb C}$.
Therefore $S_{(X, Y)}$ is not an isolated singularity.
\hfill $\Box$

\begin{remark}\label{remark:rank<=2}\rm
As in the proof of Lemma \ref{lemma:non-isolated}, we can prove that
if ${\rm rank}X \le 2$ then $S_{(X, Y)}$ is not an isolated surface singularity.
\end{remark}

\bigskip 

From Lemma \ref{lemma:non-isolated},  
we may consider $S_{(X, Y)}$ with $X \in {\rm GL}_4({\mathbb C})$ by 
changing $(X, Y)$ into $(aX+bY, cX+dY)$ for suitable 
$\left( 
\begin{array}{cc} 
a & b \\
c & d \\
\end{array}
\right) \in {\rm GL}_2({\mathbb C})$ 
for obtaining an isolated surface singularity. 
Then by choosing $P \in {\rm GL}_4({\mathbb C})$ such that ${}^{t}\!PXP = I_4$, 
we have $(X, Y) \sim ({}^{t}\!PXP, {}^{t}\!PYP) = (I_4, {}^{t}\!PYP)$. 
Hence we only need to consider the case $S_{(I_4, A)}$ with $A \in {\rm Sym}_4({\mathbb C})$.

\begin{lemma}\label{lemma:singpt}
Let $A \in {\rm Sym}_4({\mathbb C})$.  Suppose that $\dim S_{(I_4, A)} = 2$.  Then
\[{\rm Sing}(S_{(I_4, A)}) =  \{
v \in {\mathbb C}^4 \mid  {}^t v v = 0 \mbox{ and }
\exists \alpha \in {\mathbb C} \mbox{ such that } Av = \alpha v \}. \]
\end{lemma}

\noindent {\it Proof}.
The Jacobian matrix $J$ at $v \in {\mathbb C}^4$ is
$J = 2 \left(
\begin{array}{c}
{}^t v \\
{}^t (Av) \\
\end{array}
\right)$.
Hence $v \in {\mathbb C}^4$ is contained in ${\rm Sing}(S_{(I_4, A)})$
if and only if ${}^t v v = 0$ and there exists $\alpha \in
{\mathbb C} \mbox{ such that } Av = \alpha v$.
\hfill $\Box$

\begin{lemma}\label{lemma:mulroot}
Let $(X, Y) \in {\rm Sym}_4({\mathbb C}) \times {\rm Sym}_4({\mathbb C})$.
Assume that $X \in {\rm GL}_4({\mathbb C})$.
Then $(S_{(X, Y)}, O)$ is an isolated surface singularity if and only if
the characteristic polynomial of $X^{-1}Y$ has no multiple root.
Moreover, $(S_{(X, Y)}, O)$ is an isolated surface singularity if 
and only if it is a $\tilde{D}_5$-singularity. 
\end{lemma}

\noindent {\it Proof}.
For the proof, we may assume that $X=I_4$ by considering $({}^t P X P, {}^t P Y P)$ with suitable $P \in
{\rm GL}_4({\mathbb C})$. First suppose that the
characteristic polynomial of $Y$ has a multiple root $\alpha$.
By considering $Y - \alpha I_4$ instead of $Y$,
we can assume that the characteristic polynomial of $Y$ has a multiple root $0$.
If ${\rm rank} Y \le 2$, then $S_{(I_4, Y)}$ is not an isolated surface singularity
by Remark \ref{remark:rank<=2}. So suppose that ${\rm rank} \: Y=3$.
There exists $P \in {\rm GL}_4({\mathbb C})$ such that
$P^{-1}YP=J$ with
\[
J = \left(
\begin{array}{cccc}
0 & 1 & 0 & 0 \\
0 & 0 & \ast & 0 \\
0 & 0 & \ast & \ast \\
0 & 0 & 0 & \ast \\
\end{array}
\right).
\]
Set $e_1 = {}^t (1, 0, 0, 0)$ and $v = Pe_1$.
We claim that $\lambda v \in {\rm Sing}(S_{(I_4, Y)})$ for any $\lambda \in {\mathbb C}$.
From this claim we see that $S_{(I_4, Y)}$ is not an isolated surface singularity.
By Lemma \ref{lemma:singpt} it suffices to show that ${}^t v v=0$ and that
$Y v = 0$.
Since $P^{-1}YP=J$, we have $Yv=YPe_1=PJe_1 = P0=0$.
Because $Y = PJP^{-1}$ is symmetric, ${}^{t}\!(P^{-1}){}^{t}\!J {}^{t}\!P = {}^t Y = Y = PJP^{-1}$.
Hence ${}^{t}\!J {}^{t}\!P P = {}^{t}\!P P J$.
Note that $v = Pe_1$ is the second column vector of $PJ$.
So we have ${}^t v v = {}^{t}\!(Pe_1) v = {}^t e_1 {}^{t}\!P \:( \mbox{the second column of } PJ)
= {}^t e_1 \:( \mbox{the second column of } {}^{t}\!P P J) = {}^t e_1 \:( \mbox{the second column of } {}^{t}\!J {}^{t}\!P P) =0$.

Next suppose that the characteristic polynomial of $Y$ has no multiple root.
For proving that $S_{(I_4, Y)}$ is an isolated surface singularity, it suffices to prove that
if ${}^t v v=0$ for $v \in {\mathbb C}^4 \setminus \{ 0 \}$ then $v$ is not an eigenvector of $Y$
(consider the Jacobian matrix $J$ at $v$ given in the proof of Lemma \ref{lemma:singpt}).
Under the assumption, suppose that $Yv = \mu v$.
Change $Y$ into $Y-\mu I_4$. Then $Y v =0$.
There exists $P \in {\rm GL}_4({\mathbb C})$ such that
$PYP^{-1}=J$ with
\[
J = \left(
\begin{array}{cccc}
0 & 0 & 0 & 0 \\
0 & \alpha & 0 & 0 \\
0 & 0 & \beta & 0 \\
0 & 0 & 0 & \gamma \\
\end{array}
\right),
\]
where $\alpha, \beta, \gamma$ are distinct and non-zero.
Since $0 = Y v = P^{-1}JPv$, $Pv$ must be $\delta e_1$ with some $\delta \neq 0$.
Because $Y = P^{-1}JP$ is symmetric, $P{}^{t}\!PJ=JP{}^{t}\!P$
and $0 = P{}^{t}\!PJ e_1 =JP{}^{t}\!P e_1$. So we have $P{}^{t}\!P e_1 = \varepsilon e_1$ with some
$\varepsilon \neq 0$.
Hence we see that ${}^t v v = {}^{t}\!(P^{-1} \delta e_1) (P^{-1} \delta e_1) =
\delta^2 \; {}^t e_1 (P{}^{t}\!P)^{-1} e_1 = \delta^2 \varepsilon^{-1} \;  {}^t e_1 e_1 \neq 0$.
This is a contradiction. Therefore $S_{(I_4, Y)}$ is an isolated surface singularity.

Let us show that if $S_{(I_4, Y)}$ is an isolated surface singularity then it is a 
$\tilde{D}_5$-singularity. 
Note that $S_{(I_4, Y)} \subset {\mathbb C}^4$. Take a blow up of
${\mathbb C}^4$ at the origin. Let us consider the strict transform
$\tilde{S}$ of $S_{(I_4, Y)}$. We can easily check that $\tilde{S}$
is non-singular in the same way as the discussion above. The
exceptional curve $E$ is an elliptic curve defined by the two
quadratic associated with $(I_4, Y)$ in ${\mathbb P}^3$. We also see
that $E^2 = -4$. Thus $S_{(I_4, Y)}$ is a simple elliptic
singularity of $\tilde{D}_5$. \hfill $\Box$

\bigskip 

Let us describe the set $\{ (X, Y) \mid S_{(X, Y)} \text{ is a $\tilde{D}_5$-singularity } \}$ 
and the set $\{ A \mid S_{(I_4, A)} \text{ is a $\tilde{D}_5$-singularity } \}$.    

\begin{defn}\rm 
For a $4\times 4$ matrix $X$, we denote the characteristic polynomial
of $X$ by 
\[ 
P(X, t) = t^4 - c_1(X) t^3+c_2(X) t^2 -c_3(X) t +c_4(X).
\] 
Note that $c_1(X) = {\rm tr}(X)$ and $c_4(X) = \det(X)$.
Putting $a=c_1(X), b=c_2(X), c=c_3(X), d=d(X)$ for simplicity, we obtain the 
discriminant $D(X)$ of the characteristic polynomial $P(X, t)$ as 
\begin{eqnarray*} 
D(X) & = & -27a^4d^2+18a^3bcd-4a^3c^3-4a^2b^3d+a^2b^2c^2+144a^2bd^2 \\ 
& & -6a^2c^2d-80ab^2cd +18abc^3+16b^4d-4b^3c^2-192acd^2 \\
& & -128b^2d^2+144bc^2d-27c^4+256d^3. 
\end{eqnarray*} 
\end{defn}

\begin{defn}\rm
For $4\times 4$ matrices $X, Y$, we have 
\[
\begin{array}{l}
\det(sX+tY) = \\ 
c_{4, 0}(X, Y)s^4+ c_{3, 1}(X, Y)s^3t+ c_{2, 2}(X, Y)s^2t^2+c_{1, 3}(X, Y)st^3+c_{0, 4}(X, Y)t^4, 
\end{array} 
\]
where 
\[
\begin{array}{ccl}
c_{4, 0}(X, Y) & = & \det X, \\
c_{3, 1}(X, Y) & = & c_3(X){\rm tr}(Y) - c_2(X){\rm tr}(XY) + {\rm tr}(X^2Y){\rm tr}(X) - {\rm tr}(X^3Y), \\
c_{2, 2}(X, Y) & = & c_2(X)c_2(Y)+c_2(XY) -{\rm tr}(X^2Y^2) + {\rm tr}(X^2Y){\rm tr}(Y) \\
& & + {\rm tr}(XY^2){\rm tr}(X) -{\rm tr}(XY){\rm tr}(X){\rm tr}(Y), \\
c_{1, 3}(X, Y) & = & c_3(Y){\rm tr}(X) - c_2(Y){\rm tr}(XY) + {\rm tr}(XY^2){\rm tr}(Y) - {\rm tr}(XY^3), \\
c_{0, 4}(X, Y) & = & \det Y.  \\
\end{array}
\]
Then we define $D(X, Y)$ as the discriminant of the polynomial 
$\det(tX+Y)$: 
\[
\begin{array}{l}
D(X, Y) = \\
 256 (\det X)^3 (\det Y)^3 -27(\det X)^2 c_{1, 3}(X, Y)^4 -27 (\det Y)^2 c_{3, 1}(X, Y)^4 \\
-4c_{3, 1}(X, Y)^3 c_{1, 3}(X, Y)^3 -128 (\det X)^2 (\det Y)^2 c_{2, 2}(X, Y)^2 \\
+16 \det X \det Y c_{2, 2}(X, Y)^4 -80 \det X \det Y c_{3, 1}(X, Y) c_{2, 2}(X, Y)^2 c_{1, 3}(X, Y) \\
-4 \det X c_{2, 2}(X, Y)^3 c_{1, 3}(X, Y)^2 -4 \det Y c_{3, 1}(X, Y)^2 c_{2, 2}(X, Y)^3 \\
+18 \det Y c_{3, 1}(X, Y)^3 c_{2, 2}(X, Y) c_{1, 3}(X, Y) \\ 
+18 \det X c_{3, 1}(X, Y) c_{2, 2}(X, Y) c_{1, 3}(X, Y)^3 \\
+144 \det X (\det Y)^2 c_{3, 1}(X, Y)^2 c_{2, 2}(X, Y) \\ 
+144 (\det X)^2 \det Y c_{2, 2}(X, Y) c_{1, 3}(X, Y)^2 \\
+ c_{3, 1}(X, Y)^2 c_{2, 2}(X, Y)^2 c_{1, 3}(X, Y)^2 - 6 \det X \det Y c_{3, 1}(X, Y)^2 c_{1, 3}(X, Y)^2 \\
-192 (\det X)^2 (\det Y)^2 c_{3, 1}(X, Y) c_{1, 3}(X, Y). \\
\end{array}
\]
Note that
\[ D(X, Y) = (\det X)^6 D(-X^{-1} Y) = (\det X)^6 D(X^{-1} Y) \]
holds for $X \in {\rm GL}(4, {\mathbb C})$ and $Y \in {\rm M}_4({\mathbb C})$.
\end{defn}


\begin{prop}\label{prop:dxy}
Set 
\[ 
({\rm Sym}_4({\mathbb C}) \times {\rm Sym}_4({\mathbb C}))^0 :=
\{ (X, Y) \in {\rm Sym}_4({\mathbb C}) \times {\rm Sym}_4({\mathbb C}) \mid D(X, Y) \neq 0 \} 
\] 
and 
\[
{\rm Sym}_4({\mathbb C})^{0} :=
\{ X \in {\rm Sym}_4({\mathbb C}) \mid D(X) \neq 0 \}. 
\] 
Then we have the following equalities 
\[
\begin{array}{l} 
({\rm Sym}_4({\mathbb C}) \times {\rm Sym}_4({\mathbb C}))^0   \\ 
= \{ (X, Y) \in {\rm Sym}_4({\mathbb C}) \times {\rm Sym}_4({\mathbb C}) 
\mid S_{(X, Y)} \text{ is a $\tilde{D}_5$-singularity } \}
\end{array} 
\] 
and 
\[
{\rm Sym}_4({\mathbb C})^{0} = 
\{ X \in {\rm Sym}_4({\mathbb C}) \mid S_{(I_4, X)} \text{ is a $\tilde{D}_5$-singularity } \}. 
\]    
\end{prop} 

\noindent {\it Proof}.
By Lemma \ref{lemma:mulroot} we see that $S_{(I_4, X)}$ is a $\tilde{D}_5$-singularity 
if and only if $D(X) \neq 0$. So let us show that $S_{(X, Y)}$ is a $\tilde{D}_5$-singularity 
if and only if $D(X, Y) \neq 0$.
Note that $D(X, Y)=D(Y, X)$ and that $D(\lambda X, \mu Y)=\lambda^{12}\mu^{12} D(X, Y)$ for 
$\lambda, \mu \in {\mathbb C}$. 
Since the discriminant of $\det((t+a)X+Y)$ is equal to the one of $\det(tX+Y)$, 
we have $D(X, aX+Y)=D(X, Y)$ for each $a \in {\mathbb C}$.  
Hence we easily see that $D(aX+bY, cX+dY)=(ad-bc)^{12} D(X, Y)$ for 
each $\left( 
\begin{array}{cc} 
a & b \\
c & d \\
\end{array}
\right) \in {\rm GL}_2({\mathbb C})$.  

Suppose that $S_{(X, Y)}$ is a $\tilde{D}_5$-singularity. 
Then there exist $a, b \in {\mathbb C}$ such that 
$aX+bY \in {\rm GL}_4({\mathbb C})$ by Lemma \ref{lemma:non-isolated}. 
Let us choose $c, d \in {\mathbb C}$ such that 
$\left( 
\begin{array}{cc} 
a & b \\
c & d \\
\end{array}
\right) \in {\rm GL}_2({\mathbb C})$.  
Since $S_{(X, Y)} \cong S_{(aX+bY, cX+dY)}$, the characteristic 
polynomial of $(aX+bY)^{-1}(cX+dY)$ has no multiple root by Lemma \ref{lemma:mulroot}. 
Hence $D(X, Y) = (ad-bc)^{-12}D(aX+bY, cX+dY) = (ad-bc)^{-12} \det(aX+bY)^6 D((aX+bY)^{-1}(cX+dY))
\neq 0$. 

Conversely, suppose that $S_{(X, Y)}$ is not a $\tilde{D}_5$-singularity. 
If $aX+bY$ is a singular matrix for each $a, b \in {\mathbb C}$, then 
$\det(tX+Y)=0$ in ${\mathbb C}[t]$ and hence $D(X, Y)=0$. 
So we may assume that $aX+bY \in {\rm GL}_4({\mathbb C})$ for some 
$a, b \in {\mathbb C}$. 
Because $S_{(X, Y)} \cong S_{(aX+bY, cX+dY)}$ for suitable $c, d \in {\mathbb C}$ 
and $S_{(X, Y)}$ is not a $\tilde{D}_5$-singularity,  
the characteristic 
polynomial of $(aX+bY)^{-1}(cX+dY)$ has a multiple root by Lemma \ref{lemma:mulroot}. 
Hence $D(X, Y) = (ad-bc)^{-12} \det(aX+bY)^6 D((aX+bY)^{-1}(cX+dY)) = 0$.
\hfill $\Box$

\begin{defn}\rm
We define $j$-function on ${\rm Sym}_4({\mathbb C})^0$ as follows:
For $X \in {\rm Sym}_4({\mathbb C})^0$, the pair $(I_4, X)$ defines
a $\tilde{D}_5$-singularity $(S_{(I_4, X)}, O)$ because of Proposition 
\ref{prop:dxy}. Then we have an exceptional curve $E_{(I_4, X)}$ in the
minimal resolution of $(S_{(I_4, X)}, O)$. The curve $E_{(I_4, X)}$ is an
elliptic curve, and we define $j(X) := j(E_{(I_4, X)})$.
\end{defn}

\begin{theorem}\label{th:jinvariant}
We denote by $t^4-at^3+bt^2-ct+d$ the characteristic polynomial of
$X \in {\rm Sym}_4({\mathbb C})^0$. The $j$-function $j : {\rm Sym}_4({\mathbb C})^0
\to {\mathbb C}$ is given by
\[ \displaystyle
\begin{array}{ccc}
j(X)  & = &
2^8 (1728d^3-1296acd^2+432b^2d^2+324a^2c^2d-216ab^2cd \\
& & +36b^4d-27a^3c^3+27a^2b^2c^2-9ab^4c+b^6)/D,
\end{array}
\]
where the discriminant $D$ of the characteristic polynomial is given by 
$D = -27a^4d^2+18a^3bcd-4a^3c^3-4a^2b^3d+a^2b^2c^2+144a^2bd^2 
-6a^2c^2d-80ab^2cd +18abc^3+16b^4d-4b^3c^2-192acd^2 
 -128b^2d^2+144bc^2d-27c^4+256d^3$.
\end{theorem}

\noindent {\it Proof}. 
In \cite[\S 3.3.4]{Ono}, we have 
\[
j(C(M, N)) = \frac{2^8(N^2-MN+M^2)^3}{M^2N^2(N-M)^2}, 
\]
where 
$C(M, N) = \{ (X, Y, Z) \in {\mathbb P}^2 \mid 
MX^2Y-NXY^2+(X-Y)Z^2 = 0 \}$ with 
$M \neq 0, N \neq 0, M \neq N$. 
Since 
$E(M, N) = \{ 
(x_0, x_1, x_2, x_3) \in {\mathbb P}^3 \mid 
x_0^2+Mx_1^2=x_2^2, x_0^2+Nx_1^2=x_3^3 \}$ is 
isomorphic to $C(M, N)$, 
$j(E(M, N)) = j(C(M, N))$. 

Let $\lambda_0, \lambda_1, \lambda_2, \lambda_3$ be 
eigenvalues of $X \in {\rm Sym}_4({\mathbb C})^{0}$.  
The elliptic curve $E_{(I_4, X)}$ is isomorphic to 
$\{ (x_0, x_1, x_2, x_3) \in {\mathbb P}^3 \mid 
x_0^2+x_1^2+x_2^2+x_3^2 = \lambda_0 x_0^2+ 
\lambda_1 x_1^2+ \lambda_2 x_2^2+
\lambda_3 x_3^2 = 0 \}$.
By \cite[Theorem 2.5]{Ono}, setting 
\[
M = \frac{\lambda_1 - \lambda_3}{\lambda_0 - \lambda_3} \text{ and }  
N = \frac{\lambda_1 - \lambda_2}{\lambda_0 - \lambda_2}, 
\]
we obtain $E_{(I_4, X)} \cong E(M, N)$. 
Hence we can calculate the $j$-invariant of $E_{(I_4, X)}$. 
\hfill $\Box$


\begin{thebibliography}{99}

\bibitem{Artin}
M. Artin. Lectures on deformations of singularities. 
Lectures on Mathematics and Physics 54, Tata Institute of Fundamental Research, 1976. 

\bibitem{DV}
P. Du Val. On isolated singularities which do not affect the conditions of adjunction, I,II,III, 
Proc. Cambridge Philos. Soc. 30 (1934), 453-491.

\bibitem{CH}
C. Chevalley, Invariants of finite groups generated by reflections, Amer. J. Math. 77 (1955), 778-782.

\bibitem{BW}
A. Bj\"orner and V. Welker. 
The homology of ``$k$-equal'' manifolds and related partition lattices. 
Adv. Math. 110, No.2, 277-313, 1995. 

\bibitem{Brieskorn1}
E. Brieskorn. Singular elements of semi-simple algebraic groups,
Actes Congr\`es Int. Math., 2, 279-284, 1970.

\bibitem{CLO}
D. A. Cox, J. Little, and D. O'Shea. 
Using algebraic geometry. 2nd ed. 
Graduate Texts in Mathematics 185, Springer, 2005. 

\bibitem{GHC}
V. Ginzburg, Harish-Chandra bimodules for quantized Slodowy slices, 
Represent. Theory 13 (2009), 236-271.

\bibitem{grothen}
A. Grothendieck. Seminaire C. Chevally: Anneaux de Chow et
Applications, Sec. Math\'{e}matique, Expos\'{e} No.5, Paris, 1965. 

\bibitem{HS}
S. Helmke and P. Slodowy. 
Loop groups, elliptic singularities and principal bundles over 
elliptic curves. 
Geometry and topology of caustics - CAUSTICS '02. 
Banach Cent. Publ. 62, 87-99, 2004. 

\bibitem{Looijenga2}
E. J. N. Looijenga, On the semi-universal deformations of a simple elliptic singularity II, 
Topology 17 (1978), 23-40.

\bibitem{Looijenga3}
E. J. N. Looijenga, Invariant theory for generalized root systems, 
Invent. Math. 38 (1980), 1-32.

\bibitem{Mer}
J. Y. M\'erindol, Les singularit\'es simples elliptiques, leurs d\'eformations, les surfaces de Del
Pezzo et les transformations quadratiques, Ann. Scient. Ecole Norm. Sup. (4) 15 (1982), 17Ð44.

\bibitem{CRAS}
K. Nakamoto and M. Tosun, 
Semi-universal deformation spaces of some simple elliptic singularities,  
C. R., Math., Acad. Sci. Paris 345, No. 1, 31-34, 2007. 

\bibitem{Ono}
T. Ono. Variations on a Theme of Euler:  
Quadratic forms, elliptic curves, and Hopf maps, 
The university series in mathematics, 
Plenum Publishing Corporation, 1994.  

\bibitem{Pink}
H. C. Pinkham, Simple elliptic singularities, Del Pezzo surfaces and Cremona transformations, 
Proc. Symp. Pure Math. 30 (1977), 69Ð71.

\bibitem{P}
A. Premet, Enveloping algebras of Sldodwy slices and the Joseph ideal, 
J. Eur. Math. Soc. 9 (2007), no.3, 487-543.

\bibitem{Saito2}
K. Saito. Einfach-elliptische Singularit\"{a}ten.
Invent. Math. 23, 289-325, 1974. 

\bibitem{Saito3}
K. Saito and D. Yoshii. Extended affine root system IV (Simply laced Elliptic Lie Algebras),  Publ. RIMS, Kyoto Univ. 36 (2000), 385-421.

\bibitem{Slodowy1}
P. Slodowy. Simple singularities and simple algebraic group, Lecture Notes in Math. 815 (1980)

\bibitem{Slodowy2}
P. Slodowy. Singularitaeten, Kac-Moody-Lie algebren, assozierte Gruppen und Verallgemeinerungen, Habilitationsschrift, 1984.

\bibitem{Tjurina}
G. N. Tjurina. Locally semiuniversal flat deformations of isolated
singularities of complex spaces, Izv. Akad. Nauk SSSR, Ser. Mat.
33, No. 5, 967-999, 1970



\end{thebibliography}
\end{document}